\documentclass[11pt]{amsart}

\usepackage{amssymb,amsmath,latexsym, amscd, graphicx}
\usepackage{epstopdf}
\newtheorem{theorem}{Theorem}[section]
\newtheorem{lemma}[theorem]{Lemma}
\newtheorem{proposition}[theorem]{Proposition}

\def\b{\beta}
\def\e{\epsilon}
\def\s{\sigma}
\def\Q{\mathbb{Q}}
\def\Z{\mathbb{Z}}

\begin{document}

\title{ Densely ordered braid subgroups}

\date{\today}

\author[Adam Clay]{Adam Clay}
\address{Department of Mathematics\\
University of British Columbia \\
Vancouver \\
BC Canada V6T 1Z2}
\email{aclay@math.ubc.ca}
\urladdr{http://www.math.ubc.ca/~aclay/}

\author[Dale Rolfsen]{Dale Rolfsen}
\address{Department of Mathematics\\
University of British Columbia \\
Vancouver \\
BC Canada V6T 1Z2}
\email{rolfsen@math.ubc.ca}
\urladdr{http://www.math.ubc.ca/~rolfsen/}

\begin{abstract}
Dehornoy showed that the Artin braid groups $B_n$ are left-orderable.  This ordering is discrete, but we show that, for $n >2$ the Dehornoy ordering, when restricted to certain natural subgroups, becomes a dense ordering.  Among subgroups which arise are the commutator subgroup and the kernel of the Burau representation (for those $n$ for which the kernel is nontrivial).  These results follow from a characterization of least positive elements of any normal subgroup of $B_n$ which is discretely ordered by the Dehornoy ordering.
\end{abstract}

\maketitle
 
\section{Left-orderable groups and braids}

A group $G$ is \emph{left-orderable} if there exists a strict total ordering $<$ of its elements such that $g<h$ implies $fg < fh$ for all $f,g,h \in G$.  If $(G,<)$ is a left-ordered group, then the positive cone
$P := \{ g \in G \colon 1<g\}$, where $1$ is the identity of $G$, satisfies:
\begin{enumerate}
\item P is closed under multiplication and
\item For each $g \in G$, exactly one of $g=1$, $g \in P$ or $g^{-1} \in P$ holds.
\end{enumerate} 
Conversely, as is well-known, the existence of a subset $P$ satisfying these conditions implies that $G$ is left-orderable, by defining $g<h$ if and only if $g^{-1}h \in P$.  Left-orderable groups (the same class as right-orderable groups) are torsion-free and moreover they obey the zero-divisor conjecture:
if $R$ is an integral domain, then the group ring $RG$ has no zero divisors.  For further information on left-orderable groups, see \cite{MuraRhemtulla77, KopytovMedvedev96}.  A left-ordering $<$ is said to be \emph{discrete} if the positive cone has a least element.  The following is routine to verify.

\begin{proposition}
If a left ordering $<$ on $G$ is discrete, with least positive element $\e$, then every element $g \in G$ has the immediate predecessor
$g\e^{-1}$ and immediate successor $g\e$.  If it is not discrete, then it is dense, in the sense that whenever $f < g$ in $G$, there exists $h \in G$ with $f < h < g$.
\end{proposition}

For each integer $n \ge 2$, the Artin braid group $B_n$ is the group generated by 
$\s_1 , \s_2 , \dots , \s_{n-1}$, subject to the relations 
$$\s_i\s_j = \s_j\s_i {\rm \: if \:} |i-j| >1,\quad \s_i\s_j\s_i=\s_j\s_i\s_j {\rm \: if \: } |i-j| =1.$$  

It was shown by Dehornoy (see \cite{Dehornoy94, DDRW02}) that each $B_n$ is left-orderable.  The positive cone consists of all elements expressible as a word in the $\s_i$ such that the generator with the greatest subscript occurs with only positive exponents.   More precisely, for $i \in \{1, \dots n-1\}$, a word in the generators of $B_n$ is called $i$-positive (respectively $i$-negative) if it contains only the generators $ \s_1, \dots, \s_i $, the generator $\s_i$ occurs in the word, and every occurence of $\s_i$ has positive (resp. negative) exponent.  If the word contains no $\s_j, \; j \ge i$ (i.e., contains only $\s_1 \dots, \s_{i-1})$, it is said to be $i$-neutral.  A braid $\b \in B_n$ is then said to be $i$-positive (resp. $i$-neutral, $i$-negative) if it admits a representative word that is $i$-positive (resp. $i$-neutral, $i$-negative). According to Dehornoy, these possiblities are mutually exclusive: for each braid group element $\b \ne 1$
and each $i \in \{1, \dots n-1\}$, if $\b \in \langle \s_1, \dots, \s_i\rangle$ there is exactly one of three possibilities: $\b$ is $i$-positive, $i$-negative or $i$-neutral.  The positive cone therefore consists of all braids which are $i$-positive for some $i$.

\begin{proposition}
The Dehornoy ordering of $B_n$ is discrete, with smallest positive element $\s_1$.  
\end{proposition}

\begin{proof} 
Clearly $\s_1 > 1$.  Suppose there exists $\b \in B_n$ with $1 < \b < \s_1.$  Then 
$\b$ is $i$-positive for some $1 \le i \le n-1$.  If $i > 1$, then $\s_1^{-1}\b$ is also $i$-positive, so 
$\s_1^{-1}\b > 1$ and $\b > \s_1$, which
contradicts $\b < \s_1$.  If $i = 1$, then $\b$ must be a power of $\s_1$.  But the requirement that $\b < \s_1$ implies this power cannot be positive, again a contradiction.
\end{proof}

\section{Discretely ordered normal subgroups of $B_n$}

It may come as a surprise that a discretely left-ordered group can have a subgroup upon which the restriction of the ordering is dense.  A simple example of this phenomenon is the direct product $\Q \times \Z$ of the (additive) rationals and integers.  The product is given the lexicographical ordering, that is $(p,m) < (q,n)$ if and only if either $p<q$ or else $p=q$ and $m<n$, using the standard ordering of $\Q$ and $\Z$.  This ordering is discrete (and both left- and right-invariant), with least positive element 
$\e = (0,1)$.  On the other hand, the subgroup $\Q \times \{0\}$, under the same ordering is densely ordered.  

The main point of this article is that this phenomenon also happens for certain natural subgroups of the braid groups, with the Dehornoy ordering; it becomes dense when restricted to the subgroup.  To analyze this, it is first necessary to understand discretely ordered subgroups, and in particular, normal ones.

In this section, we characterize the possible least positive elements of normal subgroups of $B_n$ for which the Dehornoy ordering is discrete.  This will be used later to show that several important subgroups of 
$B_n$ are actually densely ordered by the Dehornoy ordering.  We recall the natural inclusions 
$B_m \subset B_n$ whenever $m \le n$ which takes $\s_i \in B_m$ to $\s_i \in B_n$.

In the following proofs, we use the notation $C(r)$ to denote the centralizer of $B_{r-1}$ in $B_{r}$. 
\begin{lemma}
\label{key}
Suppose that $N$ is a nontrivial normal subgroup of $B_n$ with $N \cap B_{n-1} = 1$, and $n \geq 3$.  Then if $N$ is discretely ordered by the Dehornoy ordering, its least positive element is contained in $C(n)$.
\end{lemma}
\begin{proof}
 Choose $\beta > 1$ in $N$.  Then $\beta$ must be $(n-1)$-positive, else it would lie in the trivial intersection $N \cap B_{n-1}$.  Suppose that $\beta \not\in C(n)$, so that we may choose $\gamma \in B_{n-1}$ (thus $\gamma$ is $(n-1)$-neutral) not commuting with $\beta$.

Note that $\beta \gamma \beta^{-1}$ cannot be $(n-1)$-neutral, for then $\beta \gamma \beta^{-1} \gamma^{-1}$ is also $(n-1)$-neutral, so that $\beta \gamma \beta^{-1} \gamma^{-1} \in N \cap B_{n-1} = 1$, contradicting our choice of $\gamma$.  We therefore consider two cases:
\begin{enumerate}
\item Case 1.  $\beta \gamma \beta^{-1}$ is $(n-1)$-positive.  Then $\beta \gamma \beta^{-1} \gamma^{-1}$ is also $(n-1)$-positive, so $\beta \gamma \beta^{-1} \gamma^{-1} >1$.  Since $\beta$ is $(n-1)$-positive, $\gamma \beta^{-1} \gamma^{-1}$ is $(n-1)$-negative, so $\gamma \beta^{-1} \gamma^{-1} <1$ and $\beta \gamma \beta^{-1} \gamma^{-1} < \beta$, and therefore $ 1 <  \beta \gamma \beta^{-1} \gamma^{-1} < \beta$.
\item Case 2.  $\beta \gamma \beta^{-1}$ is $(n-1)$-negative.  The inverse $\beta \gamma^{-1} \beta^{-1}$ is $(n-1)$-positive, and therefore so is $ \beta \gamma^{-1} \beta^{-1} \gamma $.  But $\gamma^{-1} \beta^{-1} \gamma$ is $(n-1)$-negative, and so $\beta > \beta \gamma^{-1}\beta^{-1}\gamma$, and we get $\beta > \beta \gamma^{-1}\beta^{-1}\gamma >1$.

\end{enumerate}

We conclude that any $\beta \notin C(n)$ cannot be the least positive element in $N$.
\end{proof}

To extend this result:

\begin{theorem} 
\label{th:leastelt}
Let $N \lhd B_n$ be a discretely ordered nontrivial normal subgroup.  Then the least positive element of $N$ is either a positive power of $\s_1$, or lies in $C(r)$, where $3 \leq r \leq n-2$ is the largest integer such that $N \cap B_{r-1}$ is trivial. 
\end{theorem}
\begin{proof} Suppose that $N \lhd B_n$, and that $N \cap B_{r-1}$ is not trivial, for all $r \geq 3$. Then $N \cap B_2$ is non-trivial, so that $N$ contains a positive power of $\s_1$, say $\s_1^m$, $m \geq 1$.  Since $\s_1$ is the least positive element in $B_n$, its immediate successors are
\[\s_1^2, \s_1^3, \s_1^4, \cdots \] 
so that $N$ can contain at most $m-1$ positive elements $<\s_1^m$, and its least element is a positive power of $\s_1$.

Assume then that for some $r \ge 3$, the intersection $N \cap B_{r-1}$ is trivial, while  $N \cap B_r$ is not.  Then we may apply our lemma, with $r$ replacing $n$, to conclude that the least positive element of  $N \cap B_r$ lies in $C(r)$.

It remains to establish that the elements of $C(r)$ are minimal in all of $N$.  Given any $\beta \in N$ with $\beta >1$, suppose that $\beta$ is $i$-positive for some $i > r-1$  (The lemma has already dealt with the case $i=r-1$).  Choose $\gamma>1$ in $N$ that is $r$-positive (or $j$-positive for any $i >j \geq r-1$).  Then the braid $\gamma^{-1} \beta$ is 
$i$-positive, so that $1< \gamma^{-1} \beta$, and we get $\gamma < \beta$.  This gives $1<\gamma<\beta$. 

\end{proof}
 
It is a result of \cite{FRZ96} that $C(r)$ consists of all elements of the form
\[ (\s_{2}\s_{1}^2 \s_{2})^p \s_{1}^q \hspace{1em} \mbox{if $r=3$}, \]
 and
\[ (\s_{r-1} \s_{r-2} \cdots \s_{2} \s_{1}^2 \s_{2} \cdots  \s_{r-2} \s_{r-1})^p \Delta_{r-1}^{2q} \hspace{1em} \mbox{if $r > 3$}.\]

Here, $\Delta_k$ is the Garside ``half-twist'' braid 
$$\Delta_k := (\s_{k-1} \s_{k-2} \cdots \s_{1})(\s_{k-1} \s_{k-2} \cdots \s_{2}) \cdots (\s_{k-1}\s_{k-2})(\s_{k-1})$$
whose square generates the (infinite cyclic) centre of $B_k$.  Note that in both cases above, $C(r) \cong \mathbb{Z} \times \mathbb{Z}$, as the two parts of each expression commute with one another.

We may rewrite the elements of $C(r)$ as 
\[ \Delta_3^{2u} \s_{1}^v \hspace{1em} \mbox{if $r=3$}, \]
and
\[\Delta_{r}^{2u} \Delta_{r-1}^{2v} \hspace{1em} \mbox{if $r >3$},\]
using the isomorphism $\mathbb{Z} \times \mathbb{Z} \rightarrow \mathbb{Z} \times \mathbb{Z}$ given by $u=p$ and $v=q-2p$ if $r=3$, and $u=p, v=q-p$ otherwise.  Taking $u \geq 1$ in the formulas above, and $r$ as in Theorem \ref{th:leastelt}, we have a description of all possible forms for the least positive element in a discretely ordered normal subgroup $N \lhd B_n$.  The condition $u \geq 1$ restricts us to considering only the positive braids described by these formulas.

\begin{proposition}
\label{prop:leastelt} Let $N \lhd B_n$ be a nontrivial normal subgroup, and let $r$ be as in Theorem \ref{th:leastelt}.  If $N$ has a least positive element $\beta$, it is of the form 
\[ \beta=\Delta_{r}^{2u}, \hspace{1em} \mbox{$u \geq 1$}\]
 or 
\[ \beta=\s_{1}^u, \hspace{1em}  \mbox{$u \geq 1$}.\] 


\end{proposition}
\begin{proof}
We proceed by ruling out all cases except for those listed above.


Suppose $\beta = \Delta_{r}^{2u} \Delta_{r-1}^{2v}$ with $v \ne 0$, and consider the normal closure $N^{\prime}$ of $\{\b\}$ in $ N$.  We show that $N^{\prime}$ contains a positive $\gamma$ with $\gamma < \beta$.

If $v<0$, take $\gamma = \beta \s_{r-1} \beta^{-1} \s_{r-1}^{-1}$.  Then $\s_{r-1} \beta^{-1} \s_{r-1}^{-1} <1$ can be seen by cancelling $\s_{r-1}$ with one of the copies of $\s_{r-1}^{-1}$ that appears in the formula for the inverse of $\Delta_{r}^{2u}$.  Therefore $\beta \s_{r-1} \beta^{-1} \s_{r-1}^{-1} < \beta$.

To show that $\gamma >1$, we compute
\[ \gamma =  \Delta_{r}^{2u} \Delta_{r-1}^{2v} \s_{r-1}  \Delta_{r}^{-2u} \Delta_{r-1}^{-2v} \s_{r-1}^{-1} \]
which can be simplified by cancelling powers of the element $\Delta_{r}^{2}$, which commute with the other parts of $\gamma$, to give
\[ \gamma =\Delta_{r-1}^{2v} \s_{r-1}  \Delta_{r-1}^{-2v} \s_{r-1}^{-1}.\]
Writing $\alpha : =\Delta_{r-1}^{2v}$, and drawing only the first $r$ strands of the braid $\gamma$, we see via an isotopy that $\gamma$ must have the same sign as $\alpha^{-1}$ (See Figure \ref{fig:shift}).

\begin{figure}[ht!]
\begin{displaymath}
\begin{array}{ccc}
\begin{picture}(60, 110)
\put(10, 0){\reflectbox{\includegraphics[scale=0.25]{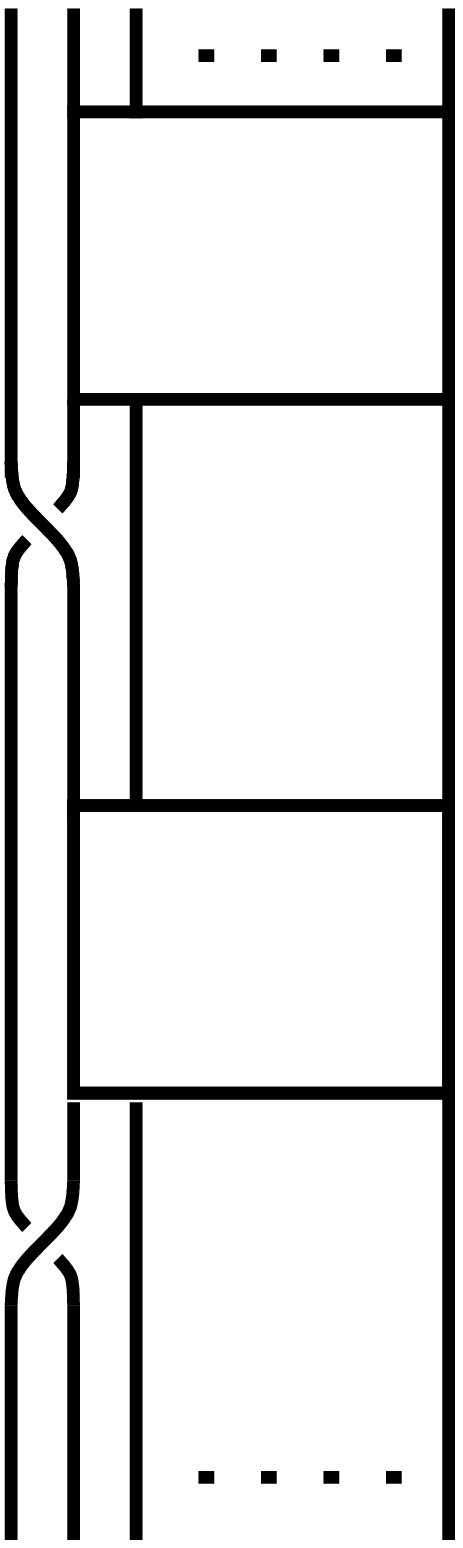}}}
\put(22, 92){$\alpha$}
\put(16, 40){$\alpha^{-1}$}
\end{picture}
 \hspace{5em} & \begin{picture}(60, 110)
\put(10, 0){\reflectbox{\includegraphics[scale=0.25]{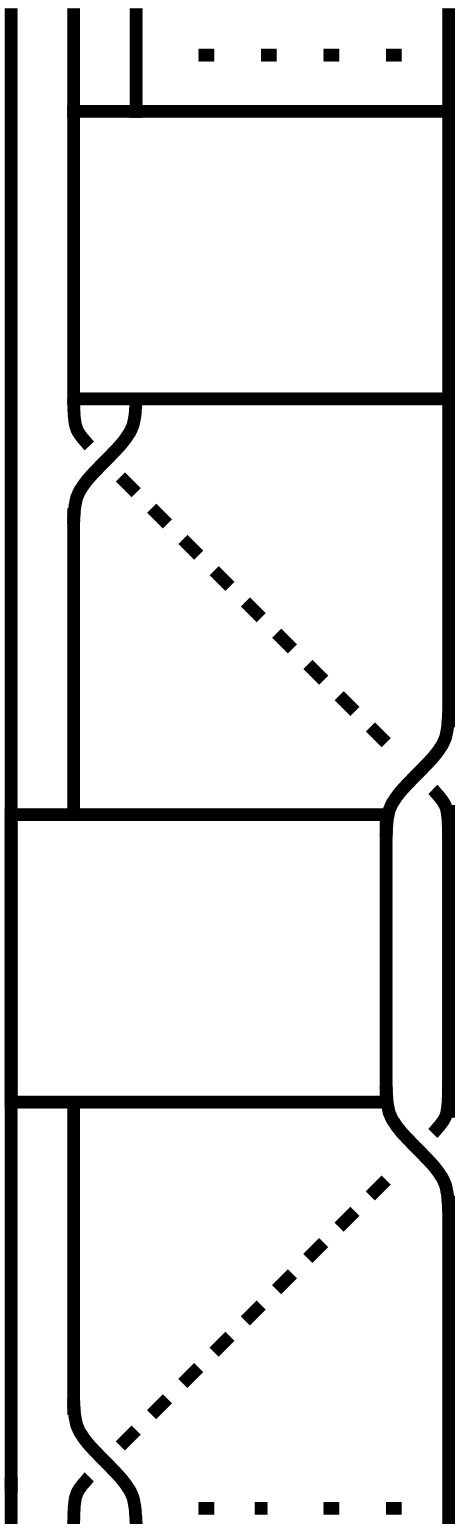}}}
\put(22, 91){$\alpha$}
\put(19, 39){$\alpha^{-1}$}
\end{picture} \\
\end{array}
\end{displaymath}
\caption{The appearance of the commutator $\alpha \s_{r-1} \alpha^{-1} \s_{r-1}^{-1}$, before isotopy on the left, and after isotopy on the right.}
\label{fig:shift}
\end{figure}

We conclude that $\gamma$ is positive if $v<0$.

If $v>0$, take $\gamma = \beta^{-1} \s_{r-1} \beta \s_{r-1}^{-1}$.  Then $\beta^{-2} \s_{r-1} \beta \s_{r-1}^{-1}$ is $(r-1)$-negative, which we can see by writing \[\beta^{-2} \s_{r-1} \beta \s_{r-1}^{-1}=\Delta_{r}^{-4u} \Delta_{r-1}^{-4v} \s_{r-1} \Delta_{r}^{2u} \Delta_{r-1}^{2v} \s_{r-1}^{-1}\]
and cancelling to get
\[\beta^{-2} \s_r \beta \s_r^{-1}= \Delta_{r-1}^{-4v}\s_{r-1} \Delta_{r}^{-2u} \Delta_{r-1}^{2v} \s_{r-1}^{-1}.\]
In the above braid we may cancel $\s_{r-1}$ with a $\s_{r-1}^{-1}$, and since $\s_{r-1}^{-1}$ appears at least twice in the formula for $\Delta_{r}^{-2u}$ (with no positive exponent occurences), the resulting braid is $(r-1)$-negative.  We conclude that $\gamma <\beta$.

To show that our choice of $\gamma$ is positive, we make an identical pictorial argument, concluding that $\gamma$ shares the same sign as $\Delta_{r-1}^{2v}$, so $\gamma$ is positive if $v>0$.

To deal with the exceptional case 
\[ \beta = \Delta_3^{2u} \s_{1}^v \]
we may use arguments identical to the two cases above, with the obvious small modifications.

\end{proof}

\section{Examples}

\subsection{The commutator subgroup}

\begin{proposition} The commutator subgroup $[B_n, B_n]$ is densely ordered by the Dehornoy ordering if $n \ge 3$.
\end{proposition}
\begin{proof}
Note that for any $n$, $[B_n, B_n]$ is characterized as the set of all braids having total exponent equal to zero for some (hence every) representative braid word.  Using this characterization, it is easy to see that for any fixed $n$ the intersection $[B_n, B_n] \cap B_{r}$ is nontrivial for all $r >2$  and trivial for $r=2$, since $[B_n, B_n]$ contains no powers of $\s_{1}$. 

By Proposition \ref{prop:leastelt}, if $[B_n, B_n]$ is to have a least positive element, it must be of the form $\Delta_3^{2u}$, $u \geq1$.  This is not possible, since $\Delta_3^{2u}$ has positive exponent sum and so is not a commutator.
We conclude that the commutator subgroup must not have a smallest positive element, and is densely ordered.
\end{proof}

The previous proposition can be improved by noting that none of the possible least elements yielded by Proposition \ref{prop:leastelt} are contained in $[B_n, B_n]$, and so a nearly identical argument applies to any normal subgroup of $B_n$ contained in the commutator subgroup.

\begin{proposition}
\label{prop:inker}
If $N$ is a normal subgroup of  $B_n$, $n \ge 3$ and $\{ 1\} \ne N \subset [B_n, B_n]$, then $N$ is densely ordered by the Dehornoy order. 
\end{proposition}

We recall that the pure braid groups $P_n$ are the subgroups of $B_n$ of braids whose associated permutation is trivial.  The Dehornoy left-ordering of $P_n$ is discrete, with least positive element 
$\sigma_{1}^2$.  On the other hand, $P_n$ has a two-sided invariant ordering \cite{RZ98}, which is necessarily dense (see \cite{LRR06}).  Noting that $P_n$ and $[P_n, P_n]$ are normal in $B_n$, and moreover
$[P_n, P_n] \subset [B_n, B_n]$ we conclude from Proposition \ref{prop:inker}: 
 
\begin{proposition}
The subgroup $[P_n, P_n]$ is densely ordered by the Dehornoy order if $n \ge 3$. 
\end{proposition}

\subsection{Brunnian braids}

An $n$-strand braid is said to be {\em Brunnian} if for every strand, the result of removing that strand is the trivial $(n-1)$ braid, that is, the identity in $B_{n-1}$.  The set of Brunnian braids is a normal subgroup of $B_n$.  For $n = 2$, it is all of $B_2$, but for $n \ge 3$ one can check that the linking number each pair of strands must be zero, which implies that Brunnian braids lie in the commutator subgroup.   As a consequence of Proposition \ref{prop:inker}, we have the following.

\begin{proposition} 
For $n \ge 3$, the subgroup of Brunnian braids in $B_n$ is densely
ordered by the Dehornoy ordering.
\end{proposition}

\subsection{Homotopically trivial braids}  

In answer to a question of Artin, D. Goldsmith \cite{Goldsmith74}
showed that for $n \ge 3$ there exist nontrivial braids which are
homotopically trivial, in the sense that they can be deformed to the
trivial braid by a homotopy of the strings, in which a string may pass
through itself, but not another string.  The set of such braids forms
a subgroup of $B_n$ which is easily seen to be normal, and in fact 
Goldsmith gives a set of generators, all of which are commutators of
pure braids.  We conclude from this (or by observing that the
candidates for least positive element are homotopically nontrivial):

\begin{proposition} 
For $n \ge 3$, the subgroup of homotopically trivial braids in $B_n$ is densely
ordered by the Dehornoy ordering.
\end{proposition}

\subsection{The Burau representation}

The Burau representation \cite{Burau36} is a map $\rho_n : B_n \rightarrow GL_n(\mathbb{Z}[t, t^{-1}])$, defined on generators by
\[\s_i \mapsto I_{i-1} \oplus
\left( \begin{array}{cc}
1-t & t\\
1 & 0  \end{array} \right) \oplus I_{n-i-1},\]
where the $1-t$ entry always appears in the $(i,i)$-th position. 

\begin{proposition} The kernels of the Burau representations $\rho_n$, where nontrivial, are densely ordered by the Dehornoy order.
\end{proposition}
\begin{proof}
 Note that $\det(\rho_n(\s_i)) = -t$ for any $i <n$. Therefore, if $\rho_n(\beta) =I_n$ for some braid $\beta$, a necessary condition is that  $\det(\rho_n(\beta))=1$, which can only happen if the total exponent sum of $\beta$ is zero.  We conclude that $\ker(\rho_n) \subset [B_n, B_n]$, and apply Proposition \ref{prop:inker}.
\end{proof}

\subsection{The homomorphism $B_4 \to B_3$}

There is a well-known homomorphism $h : B_4 \rightarrow B_3$ defined by
$$h(\s_1) = \s_1, \quad h(\s_2) = \s_2, \quad h(\s_3) = \s_1.$$
The kernel of $h$ is the normal closure of $\s_1\s_3^{-1}$ in $B_4$ and therefore lies in $[B_4,B_4]$. We conclude, from Proposition \ref{prop:inker}, the following.

\begin{proposition} The kernel of the homomorphism $h : B_4 \rightarrow B_3$ defined above is  densely ordered by the Dehornoy order.
\end{proposition}

\subsection{The Shepperd subgroup}
In contrast to the above, we consider the subgroup $H_n < B_n$ of ``braids which can be plaited with their threads tied together at 
each end,'' which was introduced by Shepperd  \cite{JS61}.
For ease of exposition, we introduce the shift homomorphism $sh \colon B_m \to B_n, \; m < n$ defined by $sh(\s_i) = \s_{i+1}$.  This is clearly injective and order-preserving.  The shift may be iterated, and we note that 
$sh^r(\Delta_{n-r})^2$ generates the center of the subgroup $\langle \s_{r+1}, \dots \s_{n-1} \rangle$ of $B_n$.  Let $H_n$ be the subgroup of $B_n$ generated by the elements
\[ \beta_i = \Delta_i^2 sh^i(\Delta_{n-i})^{-2}, \hspace{1em}i=1, \cdots, n-1,\]
and $\beta_n = \Delta_n^2$.  We call $H_n$ the Shepperd subgroup (after \cite{JS61}), where it is shown that $H_n$ is normal and that
\[H_n \cong \mathbb{Z} \times F_{n-1}, \] 
where the generator of $\mathbb{Z}$ is $\beta_n$, and the free group $F_{n-1}$ 
has generators $\beta_1, \cdots ,\beta_{n-1}$.  Note that $H_n \subset P_n$, as each of the generators is a pure braid.

\begin{proposition} For $n \geq 3$ the subgroup $H_n$ is discretely ordered, with least element $\beta_{n-1} = \Delta_{n-1}^2 $. 
\end{proposition}
\begin{proof}
Clearly $\beta_{n-1} \in H_n \cap B_{n-1}$.  In fact, the proposition follows from showing that $\langle \beta_{n-1} \rangle = H_n \cap B_{n-1}$, which we will now argue.  Suppose that there is some $\alpha \in H_n \cap B_{n-1}$ that is not a power of $\beta_{n-1}$.  Then there are two cases:

Case 1: The elements $\alpha$ and $\beta_{n-1}$ are powers of some common element $\gamma \in H_n \cap B_{n-1}$.  If $n=3$ then $\beta_{n-1} = \s_1^2$, so that the only possibility is $\gamma = \s_1$, which is not allowed since $\s_1 \notin H_3$.  On the other hand, if $n >3$, it is known (see, for example,
\cite{GM03}, Theorem 4.1) that every root of the central element $\Delta_{n-1}^2$ is conjugate to a power of $\delta$ or $\varepsilon$, where
$\delta = \s_1\s_2\cdots\s_{n-2}$ and $\varepsilon = \s_1^2\s_2\cdots\s_{n-2}$.  We note that
$ \delta^{n-1} = \varepsilon^{n-2} =\Delta_{n-1}^2$, and since their corresponding permutations have orders $n-1$ and $n-2$ respectively, no smaller power of either $\delta$ or $\varepsilon$ is a pure braid.  We conclude that this case does not occur.

Case 2: The elements $\alpha$ and $\beta_{n-1}$ are not powers of a common element.  Noting that $\alpha$ and $\beta_{n-1}$ commute yields
\[  \langle \alpha, \beta \rangle = \langle \alpha \rangle \times \langle \beta_{n-1} \rangle \cong \mathbb{Z} \times \mathbb{Z} \hookrightarrow H_n. \]
Recalling that $H_n \cong \mathbb{Z} \times F_{n-1}$ with generators as above, we find that both $\alpha$ and $\beta_{n-1}$ must map into the free group $F_{n-1}$, yielding a contradiction.
  
\end{proof}

\bibliographystyle{plain}
\bibliography{dense_braids}

\end{document}